\documentclass[12pt,oneside,english]{amsart}
\usepackage[latin1]{inputenc}
\usepackage{amssymb}

\makeatletter
\newtheorem{theorem}{Theorem}
\newtheorem{lemma}{Lemma}

\newtheorem{remark}{Remark}

\newtheorem{corollary}{Corollary}

\usepackage{babel}

\makeatother
\begin{document}
\baselineskip=17pt

\title[Theorems of Landau and Hardy-Littlewood type ]{Binary additive problems: theorems of Landau and Hardy-Littlewood type}

\author{Vladimir Shevelev}
\address{Department of Mathematics \\Ben-Gurion University of the
 Negev\\Beer-Sheva 84105, Israel. e-mail:shevelev@bgu.ac.il}


\begin{abstract}
With the uniform positions, we prove theorems of Landau and Hardy-Littlewood type for Goldbach, Chen, Lemoine-Levy and other binary partitions of positive integers. We also pose some new conjectures.
 \end{abstract}

\maketitle

\section{Introduction}
Let $\mathbb{P}$ be the set of primes and $\mathbb{P}_1$ be the set of odd primes. For an even number $n\geq6,$ let $n=p+q$ with $p,q\in \mathbb{P}$ be a Goldbach partition of $n.$ Denote by $g(n)$ the number of the unordered Goldbach partitions of $n.$ Furthermore, as in \cite{11}, denote by $\mathbf{N}_2(n)$ the number of such partitions with the taking into account the order of parts. Then, evidently,
\begin{equation}\label{1}
 \mathbf{N}_2(n)=\begin{cases}2g(n),\;\;if\;\;n/2\not\in \mathbb{P},\\2g(n)-1, \;\;if\;\;n/2\in \mathbb{P}.\end{cases}
 \end{equation}
  Justifying the posing of their famous conjecture
\begin{equation}\label{2}
\mathbf{N}_2(n)\sim2C_{tw}n/\ln^2n\prod_{p|n,\enskip p\in\mathbb{P}_1}\frac {p-1} {p-2},
\end{equation}
where $C_{tw}$ is the "twin constant"
$$ C_{tw}=\prod_{p\in \mathbb{P}_1}\frac {p(p-2)} {(p-1)^2}=0,66016... ,$$
Hardy and Littlewood (\cite{11}, pp.35-36) reproved the following Landau asymptotics
\begin{equation}\label{3}
\sum_{i=1}^{n}N_2(2i)\sim n^2/(2\ln^2n).
\end{equation}
In fact, using (\ref{3}), they obtained the following result.
\begin{theorem}\label{t1}
If the limit
$$\lim_{even \enskip x\rightarrow \infty}(\frac {N_2(x)\ln^2x} {x}\prod_{p|x,\enskip p\enskip odd\enskip prime}\frac {p-2} {p-1})$$
exists, then it is equal to $2C_{tw}.$
\end{theorem}
Landau himself proved (\ref{3}) using his identity for the $\varphi$-function which he obtained in \cite{13}:
\begin{equation}\label{4}
\sum_{d|n}\mu^2(d)/\varphi(d)=n/\varphi(n),
\end{equation}
while Hardy and Littlewood proved (\ref{3}) using their tauberian theorem \cite{10}. A modern interesting information about Landau formulas (\ref{3}),(\ref{4}) one can find in \cite{8},\cite{12}. In 1998, Dusart \cite{5} obtained the following excellent estimates for the prime counting function $\pi(x):$ if $x\geq355991,$ then
\begin{equation}\label{5}
x/\ln x+x/\ln^2x<\pi(x)<x/\ln x+x/\ln^2x+2.51x/\ln^3x.
\end{equation}
He also obtained the best lower estimate for $nth$ prime $p_n$ \cite{6}, such that the modern known estimates for $p_n$ are (\cite{1}):
\begin{equation}\label{6}
n\ln n+n\ln\ln n-n<p_n<n\ln n+n\ln\ln n,\enskip n\geq6.
\end{equation}
These estimates are much stronger than those ones which could be obtained from the well known estimates of Rosser \cite{16}.\newline
\indent In this paper, using our identities \cite{17} and estimates (\ref{5})-(\ref{6}), we give a new proof of (\ref{3}) and  close results for sums of numbers of Lemoine-Levy, Chen and some other known binary partitions; we also give new theorems like Theorem \ref{t1} (Theorems 4 and 6) and propose a new short way for proof of theorems of such type.
 \section{Goldbach partitions}
Let $\pi_1(x)$ be counting function of odd primes not exceeding $x.$ In (\cite{17}, Example 1) we showed  that, for $x=2n,$
\begin{equation}\label{7}
 \sum_{i=1}^{x/2}g(2i)=\sum_{3\leq p\leq x/2,\enskip p \in \mathbb{P}}\pi_1(x-p)
-\begin{pmatrix} \pi_1(x/2)\\2\end{pmatrix}.
\end{equation}
Using (5), we find

$$\sum_{i=1}^{x/2}g(2i)=\frac 1 {\ln x}\sum_{3\leq p\leq x/2,\enskip p \in \mathbb{P}}(x-p)-\frac {x^2} {8\ln^2x}+O(\frac{x^2} {\ln^3x})$$
\begin{equation}\label{8}
=\frac {x} {\ln x}\pi_1(x/2)-\frac {x^2} {8\ln^2x}-\frac 1 {\ln x}\sum_{i=2}^{\pi(x/2)}p_i+O(\frac{x^2} {\ln^3x})
\end{equation}
$$=\frac {3x^2} {8\ln^2x}-\frac 1 {\ln x}\sum_{i=2}^{\pi(x/2)}p_i+O(\frac{x^2} {\ln^3x}).$$
Furthermore, by (\ref{5})-(\ref{6}), we have
\newpage
\begin{equation}\label{9}
\sum_{i=2}^{\pi(x/2)}p_i=\sum_{i=2}^{\pi(x/2)}i\ln i+(\ln\ln x)\pi^2(x/2)/2+O(\frac{x^2} {\ln^2x}).
\end{equation}
Finally, notice that
\begin{equation}\label{10}
\sum_{i=2}^{\pi(x/2)}i\ln i=\int_2^{\pi(x/2)}t\ln tdt + O(\pi(x/2)\ln \pi(x/2))=\frac {x^2}{8\ln x} +O(\frac{x^2} {\ln^2x}).
\end{equation}
Now from (\ref{8})-(\ref{10}) we obtain theorem Landau in the form:
\begin{theorem}\label{t2}
\begin{equation}\label{11}
\sum_{i=1}^{x/2}g(2i)=\frac {x^2}{4\ln^2 x}+O(\frac{x^2\ln\ln x} {\ln^3x}).
\end{equation}
\end{theorem}
Note that (\ref{3}) follows from (\ref{11}) and (\ref{1}). $\blacksquare$
\begin{corollary}\label{c1}
$$\sum_{p+q\leq x,\enskip p\leq q,\enskip p, q\in \mathbb{P}}\ln(p+q)=\frac {x^2}{4\ln x}+O(\frac{x^2\ln\ln x} {\ln^2x}).$$
\end{corollary}
\bfseries Proof. \enskip\mdseries
Putting $x=2n,$ we prove that
$$\sum_{p+q\leq 2n,\enskip p\leq q,\enskip p, q\in \mathbb{P}}\ln(p+q)=\frac {n^2}{\ln n}+O(\frac{n^2\ln\ln n} {\ln^2n}).$$
Denoting, for $n\geq3,$
$$G(n)=\sum_{i=3}^{n}g(2i)$$
and using the summation by parts, we have
$$\sum_{k=3}^{n} g(2k+2)\ln (2k)=G(n+1)\ln(2n+2)-G(3)\ln6-\sum_{k=3}^{n} G(n+1)\ln (1+\frac 1{k}).  $$
Because of, by Theorem \ref{t2}, $G(n)/n$ is increasing since $n\geq n_0,$ then we have
$$\sum_{k=3}^{n} G(n+1)\ln (1+\frac 1{k})\leq \sum_{k=3}^{n} G(k+1)/k\leq $$ $$\sum_{k=3}^{n} G(n+1)/n+O(1)=G(n+1)+O(1).$$
Thus, by Theorem \ref{t2},
$$\sum_{k=3}^{n} g(2k+2)\ln (2k)=n^2/\ln n+O(\frac{n^2\ln\ln n} {\ln^2n})$$
and sufficiently to note that the following identity is valid
$$\prod_{p+q\leq x,\enskip p\leq q,\enskip p,\enskip q\in P}(p+q)=\prod_{k=3}^{n}(2k)^{g(2k)}. \blacksquare$$
\newpage
\section{Granville-Lune-Riele partitions}
In 1989, A. Granville, H. te Riele and J. van de Lune \cite{9} conjectured that if, for even $n$, $p=p(n)$ means the least prime such that  $n-p$  is also prime, then
\begin{equation}\label{12}
 p(n)=O(\ln^2n\ln\ln n).
\end{equation}
In connection with this, a Golgbach partition of $n,$ for the part $p(n)$ satisfying (\ref{12}), we call a Granville-Lune-Riele partition (GLR-partition) of $n.$ Denote $\nu(n)$ the number of unordered GLR-partition of $n.$ According to (\ref{7}), we have
\begin{equation}\label{13}
 \sum_{i=1}^{x/2}\nu(2i)=\sum_{3\leq p\leq C\ln^2x\ln\ln x,\enskip p \in \mathbb{P}}\pi_1(x-p),
\end{equation}
 (together with the replacing in (\ref{7}) $g(n)$ by $\nu(n),$ the summand $-\begin{pmatrix} \pi_1(x/2)\\2\end{pmatrix}$ is dropped since it corresponds to the now non-considered partitions $2i=i+i$ with prime $i).$ As in above, now we find
 \begin{equation}\label{13}
 \sum_{i=1}^{x/2}\nu(2i)=\frac {C}{2}x\ln x+O(x).
 \end{equation}
 By GLR-conjecture (\ref{12}), for every sufficiently large $n,$\enskip $\nu(n)>0.$ On the other hand, $\nu(n)\leq C\ln^2n\ln\ln n.$ The following arguments show that the GLR-conjecture is essentially unprovable. Note that the number of GLR-partitions of $2^k,\enskip k\geq3$ is not more than $C_1k^2\ln k$ and, consequently, the total number of such partitions for $2^k\leq x$ is $O(\ln^3x\ln\ln x).$ Note, furthermore, that, the largest parts of the Goldbach partitions of $2^k$ are  in interval $(2^{k-1}, 2^k)$ and do not take part in the Goldbach partitions as the largest parts of neither previous nor following powers for $k_1<k$ or $k_2>k.$ Thus, if to remove the largest parts in every GLR-representations of $2^k, k=3,4,...,$ we obtain a sequence $\mathbb{P}^*$ of primes with the counting function  $\pi^*(x)=\pi(x)-O(\ln^3x\ln\ln x).$ This sequence in the asymptotical sense essentially is indistinguishable from the sequence $\mathbb{P}$ of all primes with help of the approximation of $\pi(x)$ by $li(x),$ since, according to well known Littlewood result \cite{14}, the remaider term in the theorem of primes could not be less than $\sqrt{x}\log\log\log x/\log x$. But, according to our construction, for sequence $\mathbb{P}^*$ we have infinitely many even numbers for which the GLR-conjecture is wrong. Thus the conjecture is essentially unprovable.
 \begin{remark}
 Moreover, using only a weak approximation of $\pi(x)$ by $x/\ln x,$ by the similar arguments one can show that it is impossible to prove the bi
 \newpage
 nary Goldbach conjecture. Indeed, consider the sequence obtained by the following process. From every interval $(2^{m-1}, 2^m),\enskip m\geq3,$ we remove primes $p$ for which $2^m-p$ is a prime. Then the sequence gives the remaining primes $($see $A152451$ in \cite{18}$).$ Note that if $ A(x)$ is the counting function of this sequence, then, because of $\sum_{3\leq n\leq \log x}2^n/n^2=O(x/\ln^2x),$ we have $A(x)=\pi(x)-O(x/\ln^2x).$ It is well known that for the approximation of $\pi(x)$ by $x/\ln x$ the remainder term is, as the best, $O(x/\ln^2x).$ Therefore, using the approximation of $\pi(x)$ by $x/\ln x,$ the constructed sequence is in the asymptotical sense essentially indistinguishable from the sequence of all primes. But, according to the construction, we have infinitely many even numbers which are not expressible by sum of two primes of this sequence. Hence, in principle, it is impossible to prove the binary Goldbach conjecture using such approximation of $\pi(x).$\newline
 (See also an interesting sequence $A156537$ in which consecutive triples of terms form ``wells" with the depth $O(2^{n/3}\ln n/n^2)).$
 \end{remark}
\section{Lemoine-Levy partitions}
Let $l(n)$,\enskip odd $n\geq1,$ denote the number of decompositions of  $n $ into unordered sums of a prime and a doubled prime (Lemoine-Levy partitions). Then, using (\cite{17}, Example 3, where $2n-1=x$), we have:
$$\sum _{i=1}^{\frac {x+1} 2}l(2i-1)=\sum_{2\leq p\leq \frac{x+1} {2},\enskip p\in P}\pi(\frac{x-p} {2})+\sum_{2\leq p\leq \frac{x+1} {4},\enskip p\in P}\pi(x-2p)$$
\begin{equation}\label{15}
-\pi(\frac{x+1} {2})\pi(\frac{x+1} {4}).
\end{equation}
Now, as in Section 2, using (\ref{5})-(\ref{6}), we show that in average $g(n)$ and $l(n)$ asymptotically coincide.
\begin{theorem}\label{3}
\begin{equation}\label{16}
\sum _{i=1}^{\frac {x+1} 2}l(2i-1)=\frac {x^2}{4\ln^2 x}+O(\frac{x^2\ln\ln x} {\ln^3x}).
\end{equation}
\end{theorem}
\begin{corollary}\label{c2}
$$\sum_{p+2q\leq x,\enskip p,q\in \mathbb{P}}\ln(p+2q)=\frac {x^2}{4\ln x}+O(\frac{x^2\ln\ln x} {\ln^2x}).$$
\end{corollary}
\bfseries Proof. \enskip\mdseries The proof of Corollary \ref{c1} is repeated, if to note that
$$ \prod_{p+2q\leq 2n-1,\enskip p, q\in \mathbb{P}}(p+2q)=\prod_{k=5}^{n}(2k-1)^{l(2k-1)}. \blacksquare$$
\newpage
Moreover, using the scheme of proof of Theorem 6 (see below), it could be proved the following result which is similar to Theorem \ref{t1}.
\begin{theorem}\label{t4}
If the limit
$$\lim_{odd\enskip x\rightarrow\infty} (\enskip\frac {l(x)\ln^2x} {x}\prod_{p|x,\enskip p\in\mathbb{P}_1}\frac {p-2} {p-1}\enskip)$$
exists, then it equals to $C_{tw}.$
\end{theorem}
This theorem is a base for a conjecture which is very close to the Hardy and Littlewood conjecture (\ref{2}): for odd $n\rightarrow\infty$
\begin{equation}\label{17}
l(n)\sim C_{tw}n/\ln^2n\prod_{p|n,\enskip p\in \mathbb{P}_1}\frac {p-1} {p-2}.
\end{equation}

 See also our comments to sequences $A152460-A152461$ in \cite{18}.

\section{Chen partitions}
Let $c(n)$ denote the number of decompositions of even $n$ into prime and prime or semiprime (\cite{19}).
In 1966, Chen \cite{2} did a very important step towards Goldbach conjecture, proving that for sufficiently large $n$
\begin{equation}\label{18}
c(n)\geq aC_{tw}(\prod_{p|n,\enskip p\in\mathbb{P}_1}\frac {p-1} {p-2})n/\ln^2 n,
\end{equation}
where $a=0.098$ (\cite{2}) with improvements $a=0.67$ (\cite{3}) and 1.0974 (\cite{4}). Ross \cite{15} gave a large simplification of Chen's proof. \newline Putting in (\cite{17}, Example 2) $2n=x$ and
$g_1(2i)+g_2(2i)=c(2i),\enskip i=1,...,x/2,$ we have:
$$\sum _{i=1}^{x/2}c(2i)=\sum_{9\leq q\leq x/2,\enskip q\in \mathbb{P}_{2,1}}\pi_1(x-q)+\sum_{3\leq p\leq x/2,\enskip p\in \mathbb{P}}\pi_{2,1}(x-p)$$
\begin{equation}\label{19}
-\pi_1(x/2)(\pi_{2,1}(x/2)+1/2(\pi_1(x/2)-1))+\pi(x/2-1)+1, \enskip x\geq4,
\end{equation}
where $\pi_1(x)$ is the counting function of odd primes and $\pi_{2,1}(x)$ is the counting function of the set $\mathbb{P}_{2,1}=\{q_n\}_{n\geq1}$ of odd semiprimes, such that (cf. \cite{19}):
\begin{equation}\label{20}
\pi_{2,1}(x)=\sum_{3\leq p\leq \sqrt{x},\enskip p\in \mathbb{P}}(\pi(x/p)-\pi(p)+1).
\end{equation}
Using (5)-(6), we find
\begin{equation}\label{21}
\sum_{3\leq p\leq \sqrt{x},\enskip p\in \mathbb{P}}\pi(x/p)=x\ln\ln x/\ln x+O(x\ln\ln x/\ln^2 x),
\end{equation}\newpage
\begin{equation}\label{22}
\sum_{3\leq p\leq \sqrt{x},\enskip p\in \mathbb{P}}\pi(p)=O(x/\ln^2 x),
\end{equation}
such that, by (\ref{20}), we have
\begin{equation}\label{23}
\pi_{2,1}(x)=x\ln\ln x/\ln x+O(x\ln\ln x/\ln^2 x).
\end{equation}
Note that, by (\ref{23}),
$$ n=\pi_{2,1}(q_n)= q_n\ln\ln q_n/\ln q_n+O(q_n\ln\ln q_n/\ln^2 q_n),$$
or
\begin{equation}\label{24}
q_n=n\ln q_n/\ln\ln q_n+O(n/\ln\ln q_n).
\end{equation}
Iterating of (\ref{24}), we find
\begin{equation}\label{25}
q_n=n\ln q_n/\ln\ln q_n+O(n/\ln\ln q_n)$$ $$= \frac {n(\ln n+\ln\ln q_n-\ln\ln\ln q_n+o(1))} {\ln(\ln n+\ln\ln q_n-\ln\ln\ln q_n+o(1))}(1+o(1))
\end{equation}
\begin{equation}\label{26}
\ln q_n=\ln n+\ln\ln n(1+o(1)).
\end{equation}
Using (\ref{25})-(\ref{26}), we obtain that
\begin{equation}\label{27}
q_n= n\ln n/\ln\ln n+O(n\ln\ln n).
\end{equation}
Using (\ref{23}),(\ref{27}), from (\ref{19}) we have
$$\sum_{i=1}^{x/2}c(2i)=\frac 1 {\ln x}\sum_{9\leq q\leq x/2,\enskip q \in \mathbb{P}_{2,1}}(x-q)+\frac {\ln\ln x} {\ln x}\sum_{3\leq p\leq x/2,\enskip p \in \mathbb{P}}(x-p)$$ $$-\frac {x^2\ln\ln x} {4\ln^2x}+O(\frac{x^2} {\ln^2x})$$
\begin{equation}\label{28}
=\frac {3} {4} \frac {x^2\ln\ln x}{\ln^2x}-\frac {1} {\ln x}\sum_{i=1}^{\pi_{2,1}(x/2)}q_i-\frac {\ln\ln x} {\ln x}\sum_{i=1}^{\pi(x/2)}p_i+O(\frac{x^2} {\ln^2x}),
\end{equation}
or, taking into account (\ref{9})-(\ref{10}),
\begin{equation}\label{29}
\sum_{i=1}^{x/2}c(2i)=\frac {5} {8} \frac {x^2\ln\ln x}{\ln^2x}-\frac {1} {\ln x}\sum_{i=1}^{\pi_{2,1}(x/2)}q_i+O(\frac{x^2} {\ln^2x}).
\end{equation}
Furthermore, using (\ref{27}) and (\ref{23}), we have
$$\sum_{i=1}^{\pi_{2,1}(x/2)}q_i=\sum_{i=1}^{\pi_{2,1}(x/2)}i\ln i/\ln\ln i +O(\sum_{i=1}^{\pi_{2,1}(x/2)}i\ln\ln i)$$
$$=\int_3^{\pi_{2,1}(x/2)}\frac {t\ln t} {\ln\ln t}dt+O(\pi^2_{2,1}(x/2)\ln \ln\pi_{2,1}(x/2))$$\newpage
\begin{equation}\label{30}
=\int_3^{\pi_{2,1}(x/2)}\frac {t\ln t} {\ln\ln t}dt+O(\frac {x^2(\ln\ln x)^3} {\ln^2x}).
\end{equation}
Finally,
$$\int_3^{\pi_{2,1}(x/2)}\frac {t\ln t} {\ln\ln t}dt =(\frac {1} {2}\pi_{2,1}^2(x/2)\ln\pi_{2,1}(x/2)-\frac{1} {4}\pi_{2,1}^2(x/2))\frac {1} {\ln\ln \pi_{2,1}(x/2)}$$
\begin{equation}\label{31}
+O(\int_3^{\pi_{2,1}(x/2)}\frac {tdt} {(\ln\ln t)^2})=\frac{x^2} 8\frac{\ln\ln x} {\ln x}+O((\frac {x\ln\ln x} {\ln x})^2).
\end{equation}
Now from (\ref{29})-(\ref{31}) we obtain the following result.
\begin{theorem}\label{t5}
\begin{equation}\label{32}
\sum _{i=1}^{\frac {x} 2}c(2i)=\frac {x^2\ln\ln x}{2\ln^2 x}+O(\frac{x^2} {\ln^2x}).
\end{equation}
\end{theorem}
Similar to Corollaries \ref{c1},\ref{c2},  in view of the identity
$$\prod_{p+q\leq 2n,\enskip  p\in \mathbb{P},\enskip q\in \mathbb{P}_{2,1}\bigcup \mathbb{P}}(p+q)=\prod_{k=3}^{n}(2k)^{c(2k)}, $$
we obtain the following formula.
\begin{corollary}\label{c3}
$$\sum_{p+q\leq x,\enskip p\in \mathbb{P},\enskip q\in \mathbb{P}_{2,1}\bigcup \mathbb{P}}\ln(p+q)=\frac {x^2\ln\ln x}{2\ln x}+O(\frac{x^2\ln\ln x} {\ln^2x}).$$
\end{corollary}
Furthermore, we show that the following theorem is valid.
\begin{theorem}\label{t6}
If limit
$$\lim_{even \enskip x\rightarrow \infty}(\frac {c(x)\ln^2x} {x\ln\ln x}\prod_{p|x,\enskip p\in \mathbb{P}_1}\frac {p-2} {p-1})$$
exists, then it equals to $2C_{tw}.$
\end{theorem}
This theorem is a base for a similar to (\ref{2}) conjecture: for even $n\rightarrow\infty$
\begin{equation}\label{33}
c(n)\sim 2C_{tw}n\ln\ln n/\ln^2n\prod_{p|n,\enskip p\in\mathbb{P}_1}\frac {p-1} {p-2}.
\end{equation}
Note that inequality (\ref{18}) is weaker than (\ref{33}). Therefore, it is nature to think that a proof of the binary Goldbach conjecture will come using a form of a weaker inequality than (\ref{2}), for instance,
$$  c(n)\geq aC_{tw}(\prod_{p|n,\enskip p\in\mathbb{P}_1}\frac {p-1} {p-2})n/\ln^3 n.$$\newpage
\section{Proof of theorem 6}
By the following way of proof, one can also reprove Theorem 1 and prove Theorem 4.
We use a known statement belonging to O. Stolz (cf. \cite{7}, pp. 67-68).
\begin{lemma}\label{L1}
Let $y_n\rightarrow\infty$ and be increasing at least for $n\geq n_0:$
$$\Delta y_n=y_n-y_{n-1}>0,\enskip n\geq n_0.$$
If, for some another variable $x_n,$ there exists limit
$$\lim_{n\rightarrow \infty}\frac{\Delta x_n} {\Delta y_n}=a, \enskip |a|\leq\infty, $$
then $\lim_{n\rightarrow \infty}\frac{ x_n} { y_n}$ exists as well and equals to $a.$
\end{lemma}
\begin{lemma}\label{L2}
If, for $a_n>0,$ we have $1/a_n=O(1)\enskip (n\rightarrow\infty)$ and, for $\varepsilon>0,$
\begin{equation}\label{34}
\sum_{i=1}^{n}a_i=O(n(\ln n)^{1-\varepsilon}),
\end{equation}
then
\begin{equation}\label{35}
\sum_{i=3}^{n}\frac {\ln\ln i} {\ln^2i}ia_i\sim \frac {\ln\ln n} {\ln^2n}\sum_{i=3}^{n}ia_i.
\end{equation}
\end{lemma}
\bfseries Proof. \enskip\mdseries Putting
\begin{equation}\label{36}
\lambda(n)=\frac{\ln^2n} {\ln\ln n},\enskip y_n=\sum_{i=3}^{n}ia_i,\enskip x_n=\lambda(n)\sum_{i=3}^{n}\frac{1}{\lambda(i)}ia_i,
\end{equation}
we should prove that $\lim_{n\rightarrow \infty}\frac{ x_n} { y_n}=1,$ or, by Lemma 1, that $\lim_{n\rightarrow \infty}\frac{\Delta x_n} {\Delta y_n}=1.$
We have
$$\lim_{n\rightarrow \infty}\frac{\Delta x_n} {\Delta y_n}=\lim_{n\rightarrow \infty}(na_n)^{-1}(\sum_{i=3}^{n}\frac{ia_i}{\lambda(i)}(\lambda(n)-\lambda(n-1))+\frac{\lambda(n-1)}{\lambda(n)}na_n)$$
\begin{equation}\label{37}
=1+\lim_{n\rightarrow \infty}(\frac{\lambda(n)-\lambda(n-1)}{na_n}\sum_{i=3}^{n}\frac{ia_i}{\lambda(i)}).
\end{equation}
Note that, by (\ref{34}),(\ref{36}), we have
$$\sum_{i=3}^{n}\frac{ia_i}{\lambda(i)}\leq C\frac {n} {\lambda(n)}\sum_{i=3}^{n}a_i=O(\frac{n^2(\ln n)^{1-\varepsilon}} {\lambda(n)}),$$
or, using condition $1/a_n=O(1),$ we find
\begin{equation}\label{38}
\frac 1 {na_n}\sum_{i=3}^{n}\frac{ia_i}{\lambda(i)}=O(\frac{n(\ln n)^{1-\varepsilon}} {\lambda(n)}).
\end{equation}
Furthermore,\newpage
$$\lambda(n)-\lambda(n-1)\sim2\ln n/n\ln\ln n=2\frac {\lambda(n)} {n\ln n}. $$
Consequently, by (\ref{38}), we have
$$\frac{\lambda(n)-\lambda(n-1)}{na_n}\sum_{i=3}^{n}\frac{ia_i}{\lambda(i)}=O(\frac{1} {\ln^\varepsilon n)}$$
and the lemma follows from (\ref{37}).$\blacksquare$\newline
Using Lemmas \ref{L1} and \ref{L2}, we complete proof of the theorem. For $x=2n,\enskip n\geq3,$ put
$$x_n=\sum _{i=3}^{n}c(2i), \enskip y_n=\sum _{i=3}^{n}2i\ln\ln 2i/\ln^22i\prod_{p|i,\enskip p\in\mathbb{P}_1}\frac {p-1} {p-2}.$$
Then, by the condition, $\lim_{n\rightarrow \infty}\frac {x_n-x_{n-1}}{y_n-y_{n-1}}$ exists and, by Lemma \ref{L1}, it is sufficient to prove that
\begin{equation}\label{39}
\lim _{n\rightarrow \infty}\frac{x_n}{y_n}=2C_{tw}.
\end{equation}
Putting
$$a_n=\prod_{p|n,\enskip p\in\mathbb{P}_1}\frac {p-1} {p-2},$$
we note that $$a_n=\prod_{p|n,\enskip p\in\mathbb{P}_1}(1+\frac{1}{p-2})=O(\ln\ln n)$$ and the conditions
of Lemma \ref{L2} are satisfied. Using Theorem \ref{t5} and Lemma \ref{L2}, we have
\begin{equation}\label{40}
\lim _{n\rightarrow \infty}\frac{x_n}{y_n}=\lim _{n\rightarrow \infty}\frac {n^2} {\sum _{i=3}^{n}i\prod_{p|i,\enskip p\in\mathbb{P}_1}(1+\frac{1}{p-2})}.
\end{equation}
Furthermore,
\begin{equation}\label{41}
\sum _{i=3}^{n}i\prod_{p|i,\enskip p\in\mathbb{P}_1}(1+\frac{1}{p-2})=\sum _{i=3}^{n}i\sum_{d|i}\frac{\eta(d)} {d^*},
\end{equation}
where $\eta(n)=1,$ if $n$ is an odd square-free number, and $\eta(n)=0,$ otherwise;
$$d^*=\prod_{p|d}(p-2)^{\alpha_p}, \enskip if \enskip d=\prod_{p|d}p^{\alpha_p}.$$
Continuing transformations of expression (\ref{41}), we have
$$\sum _{i=3}^{n}i\prod_{p|i,\enskip p\in\mathbb{P}_1}(1+\frac{1}{p-2})=\sum_{d\leq n}\frac{\eta(d)} {d^*}\sum_{l\leq \frac {n}{ d}}ld$$
 $$=\sum_{d\leq n}\frac{\eta(d)d} {d^*}\sum_{l\leq \frac {n}{ d}}=\sum_{d\leq n}\frac{\eta(d)d} {d^*}\frac{n^2} {2d^2}+O(n\sum_{d\leq n}\frac{1} {d})$$\newpage
 $$=\frac{n^2} {2}\sum _{i=3}^{n}\frac{\eta(i)} {i^*i}+O(n\ln n)$$ $$=\frac{n^2} {2}\prod_{p\geq3, \enskip p\in P}(1+\frac {1} {p(p-2)})+O(n\ln n)$$
 \begin{equation}\label{42}
 =\frac{n^2} {2}\prod_{p\geq3, \enskip p\in P}(1-\frac {1} {(p-1)^2})^{-1}+O(n\ln n)=\frac{n^2}{2C_{tw}}+O(n\ln n),
 \end{equation}
    and (39) follows from (42) and (40).$\blacksquare$

\section{Conclusive remarks}
In this paper we selected only most interesting, in our opinion, binary partitions of positive integers. General results, which were obtained in \cite{17}, allow, by a similar way, to study many another binary partitions as well. Only for example, it could be proved that if, for  $x\equiv1 \pmod4,$ $q(x)\geq0$ denote the number of unordered representations of $x$ by a sum of two squares of nonnegative integers, then
\begin{equation}\label{43}
\sum_{i=0}^{\frac{x-1} {4}}q(4i+1)=\frac{\pi} {16}x+O(\sqrt{x}).
\end{equation}
This formula follows from the exact equality (cf. \cite{17}, Example 4)
$$\sum _{i=0}^{\frac {x-1} {4}}q(4i+1)=\sum_{k=1}^{\lfloor(1+\sqrt{x})/2\rfloor}\lfloor(1+\sqrt{2x-(2k-1)^2 }\enskip)/2\rfloor$$
\begin{equation}\label{44}
-\begin{pmatrix} \lfloor(1+\sqrt{x})/2\rfloor\\2\end{pmatrix}.
\end{equation}
In addition, note that the last formula could be rewrite in the form
\begin{equation}\label{45}
\sum _{i=0}^{n}q(4i+1)=\sum_{1\leq k\leq r(n)}r(2n-k^2+k)-\begin{pmatrix} r(n)\\2\end{pmatrix}, \end{equation}  where $r(n)$ is the nearest integer to $\sqrt {n+1}$ and $n=\frac {x-1}{4}.$


\begin{thebibliography}{18}
\bibitem 1 N. \enskip E.\enskip Bach, J.\enskip Shallit,  \enskip Algorithmic Number Theory , \slshape MIT Press, \upshape \bfseries 233\mdseries \enskip(1996). ISBN 0-262-02405-5.
\bibitem 2 J. \enskip R. \enskip Chen,\enskip  On the representation of a large even integer as the sum of a prime and the product of at most two primes,  \enskip \slshape Kexue Tongbao,\upshape\enskip \bfseries 17\mdseries\enskip(1966),\enskip385--386.
\bibitem 3 J.\enskip R. \enskip Chen,\enskip  On the representation of a larger even integer as the sum of a prime and the product of at most two primes,\enskip \slshape Sci. Sinica, \upshape \enskip\bfseries 16\mdseries \enskip(1973),\enskip157--176.
\bibitem 4 Y. \enskip Cai,\enskip A remark on Chen's theorem, \enskip \slshape Acta Arithm.\upshape, \enskip \bfseries 102\mdseries, \enskip no.4 (2002),\enskip339-352.
\bibitem 5  P. \enskip Dusart,\enskip \slshape Autour de la fonction qui compte le nombre de nombres premiers, \upshape\enskip Doctoral thesis for l'Universite de Limoges (1998).\newpage
\bibitem 6  P.\enskip Dusart ,\enskip  The $k$-th prime is greater than $k(\ln k + \ln\ln k-1)$ for $k\geq2,$  \enskip \slshape Math. Comp. \upshape\enskip \bfseries 68\mdseries \enskip(1999),\enskip 411--415.
\bibitem 7  G.\enskip  M.\enskip 	Fihtengolz,\enskip  \slshape  Course of differential and integral calculus,\upshape \enskip vol. I, Moscow, Nauka, 1969 (Russian).
\bibitem 8  G.\enskip Giordano,\enskip On the irregularity of the distribution of the sums of pairs of odd primes, \enskip \slshape Int. J. of Math. and Math. Sc.\upshape  \enskip\bfseries 30,\mdseries \enskip no.6 (2002), 377--381.
\bibitem 9 A.\enskip  Granville, \enskip J.\enskip  van de Lune, and H.\enskip J.\enskip J.  te Riele, \enskip  Checking the Goldbach conjecture on a vector computer, \enskip \slshape Number Theory and Application \upshape\enskip (R.\enskip  A.\enskip  Mollin, ed.), \enskip Kluwer Acad.,\enskip 1989,\enskip  423--433.
\bibitem {10} G.\enskip  H.\enskip  Hardy , \enskip J. E. Littlewood, \enskip  Tauberian theorems concerning power series and Dirichlet's series whose coefficients are positive,\enskip \slshape Proc. London Math.\upshape,\enskip  \bfseries 13\mdseries, ser.2 \enskip (1913/14),\enskip  174--192.
\bibitem {11} G. \enskip H. \enskip Hardy and J.\enskip  E.\enskip  Littlewood, \enskip Some problems of `partitio numerorum`.\enskip III: on the expression of a number as a sum of primes,\enskip \slshape Acta Math.\upshape\enskip  \bfseries 44\mdseries \enskip(1922), \enskip 1--70.
\bibitem {12} P. \enskip Haukkanen,\enskip On generalized Landau identities,\enskip \slshape Portugal. Math.\upshape \enskip  \bfseries 52\mdseries, no.1 (1995),\enskip29--38.
\bibitem {13} E.\enskip Landau, \enskip  \"{U}ber die zahlentheoretische Funktion $\varphi(n)$ und ihre Beziehung zum Goldbachschen Satz,\enskip \slshape G\"{o}ttinger Nachrichten \upshape\enskip(1900), \enskip177--186.
\bibitem {14} J. \enskip E. \enskip Littlewood,\enskip Sur la distribution des nombres premiers, \slshape C.\enskip R. Paris\upshape ,\enskip  \bfseries 158\mdseries \enskip(1914),\enskip 1869--1872.
\bibitem {15} P.\enskip M. \enskip Ross,\enskip On Chen's theorem that each large even number has the form $(p_1+p_2)$ or $(p_1+p_2p_3)$, \enskip \slshape J. London Math. Soc.\upshape\enskip \bfseries 2\mdseries\enskip (1975), \enskip500--506.
\bibitem {16} B.\enskip	Rosser,\enskip Explicit Bounds for some Functions of Prime Numbers,\enskip\slshape Amer. J. Math.\upshape \enskip \bfseries 63\mdseries \enskip (1941), \enskip211--232.
\bibitem {17} V.\enskip  Shevelev, \enskip Binary additive problems: recursions for numbers of representations,\enskip \slshape http://www.arxiv.org/abs/\upshape 0901.3102.
\bibitem {18} N.\enskip J.\enskip A.\enskip Sloane,\enskip\slshape The On-Line Encyclopedia of Integer Sequences \upshape\enskip http://oeis.org.
\bibitem {19} E.\enskip W.\enskip Weisstein,\enskip\slshape "Semiprime" \upshape From MathWorld--A Wolfram Web Resource (http://mathworld.wolfram.com)
\end{thebibliography}
\end{document}